\documentclass{llncs}
\begin{document}
\pagenumbering{arabic}
\pagestyle{headings}

\title{
	   Permutations avoiding a pattern from $S_k$ \\
	   and at least two patterns from $S_3$} 

\author{Toufik Mansour}
\institute{Department of Mathematics \\
	   University of Haifa, Haifa, Israel 31905\\
	   tmansur@study.haifa.ac.il}
\maketitle
\begin{abstract}
In this paper, we find explicit formulas or generating functions for 
the cardinalities of the sets $S_n(T,\tau)$ of all permutations in $S_n$ 
that avoid a pattern $\tau\in S_k$ and a set $T$, $|T|\geq 2$, of 
patterns from $S_3$. The main body of the paper is divided into three 
sections corresponding to the cases $|T|=2,3$ and $|T|\geq 4$. 
As an example, in the fifth section, we obtain the complete 
classification of all cardinalities of the sets $S_n(T,\tau)$ 
for $k=4$.
\end{abstract}
\section{Introduction}

Let $[k]=\{1,\dots,k\}$ be a (totally ordered) {\em alphabet} on $k$ 
letters, and let $\alpha\in [k]^m$, $\beta\in [l]^m$ with $l\leq k$. 
We say that $\alpha$ is {\em order-isomorphic} to $\beta$ if the 
following condition holds for all $1\leq i<j\leq n$: 
$\alpha_i<\alpha_j$ if and only if $\beta_i<\beta_j$.\\

We say that $\tau\in S_n$ {\em contains} $\alpha\in S_k$ 
if there exist $1\leq i_1<\dots<i_k\leq n$ such that 
$(\tau_{i_1},\dots,\tau_{i_k})$ is order-isomorphic to 
$\alpha=(\alpha_1,\dots,\alpha_k)$, and we say that $\tau$ {\em avoids} 
$\alpha$ if $\tau$ does not contain $\alpha$. 
The set of all permutations in $S_n$ avoiding $\alpha$ is denoted 
$S_n(\alpha)$. More generally, for any finite set of permutations $T$ 
we write $S_n(T)$ to denote the set of permutations in $S_n$ avoiding all 
the permutations in $T$.
Two sets, $T_1, T_2$, are said to be {\em Wilf equivalent} (or to belong to 
the same {\em Wilf class}) if and only if $|S_n(T_1)|=|S_n(T_2)|$ for 
any $n\geq 0$; the Wilf class of $T$ we denote by $\overline{T}$.\\

The study of the sets $S_n(\alpha)$ was initiated by Knuth ~\cite{knuth}, 
who proved that $|S_n(\alpha)|=\frac{1}{n+1}{{2n}\choose n}$ for any 
$\alpha\in S_3$. Knuth's results where further extended in two directions. 
West ~\cite{west} and Stankova ~\cite{stankova} analyzed $S_n(\alpha)$ for 
$\alpha\in S_4$ and obtained the complete classification, which contains $3$ 
distinct Wilf classes. This classification, however, does not give exact values 
of $S_n(\alpha)$. On the other hand, Simion and Schmidt ~\cite{simion} 
studied $S_n(T)$ for arbitrary subsets $T\subseteq S_3$ and 
discovered $7$ Wilf classes. The study of $S_n(\alpha,\tau)$ 
for all $\alpha\in S_3$, $\tau\in S_4(\alpha)$, was 
completed by West ~\cite{west}, Billey, Jockusch and Stanley ~\cite{billey}, 
and Guibert ~\cite{guibert}. \\

In the present paper, we calculate the cardinalities of the sets 
$S_n(T,\tau)$ for all $T\subseteq S_3$, $|T|\geq 2$, and all 
permutations $\tau\in S_k$ such that $k\geq 3$.

\begin{remark}
\label{rem}
West ~\cite{west} observed that if $\tau$ contains a pattern in 
$T$, then $|S_n(T,\tau)|=|S_n(T)|$.
 Therefore, in what follows we assume that $\tau\in S_k(T)$.
\end{remark}

Throughout the paper, we often make use of the following simple
statement.

\begin{lemma}
\label{det}
Let $\{s_i(x)\}_{i=1}^r$, $\{A_i(x)\}_{i=1}^r$ and $\{B_i(x)\}_{i=1}^r$ 
be sequences of functions such that 
		$$s_i(x)=A_i(x)s_{i+1}(x)+B_i(x),$$
where $1\leq i\leq r-1$, and $s_r(x)=h(x)$. Then 
	$$s_1(x)=\left\bracevert
		\begin{array}{ccccc}
		B_1(x)		&	-A_1(x)		&	0	&	\dots	&	0 \\
		B_2(x)		&	1		&   -A_2(x)	&	\dots	&	0 \\
		B_3(x)		&	0		&	1	&  	\dots	& 	\vdots \\
		\vdots		&	\vdots		&	0	&	\ddots	&	0	\\
		B_{r-1}(x)	&	0		&	0	&	\ddots	&	-A_{r-1}(x) \\
		h(x)		&	0		&	0	&	\dots	&	1 
		\end{array}
		\right\bracevert.$$
\end{lemma}
\begin{proof}
Immediately, by definitions and induction on $r$.\hfill\rule{2mm}{2mm}
\end{proof}

Our calculation is divided into three sections corresponding to the cases 
$|T|=2,3$ and $|T|\geq 4$. In the last section, as an example, we will 
obtain the complete classification (Table \ref{tab34}) of all cardinalities 
of the sets $S_n(T,\tau)$ where $T\subseteq S_3$, $\tau\in S_4$.
\section{Avoiding a pair from $S_3$ and a pattern from $S_k$}

In this section, we calculate the cardinality of the sets 
$S_n(\beta^1,\beta^2,\tau)$ where $\beta^1,\beta^2\in S_3$, 
$\tau\in S_k$, $k\geq 3$. By Remark \ref{rem} and by the three natural operations, the complementation, the reversal 
and the inverse (see Simion and Schmidt ~\cite{simion}, Lemma $1$), 
we have to consider the following four possibilities:\\

$\begin{array}{rlll}
  1)&	S_n(123,132,\tau),	& \ \mbox{where}	& \ \tau\in S_k(123,132), \\
  2)&	S_n(123,231,\tau), 	& \ \mbox{where}	& \ \tau\in S_k(123,231), \\
  3)&	S_n(132,213,\tau),	& \ \mbox{where}	& \ \tau\in S_k(132,213), \\
  4)&	S_n(213,231,\tau),	& \ \mbox{where} 	& \ \tau\in S_k(213,231). \\
\end{array}$
\\

The main body of this section is divided into four subsections corresponding 
to the above four cases. 
\subsection{$T=\{123,\ 132\}$.}

Let $a_{\tau}(n)=|S_n(123,132,\tau)|$, and let $a_{\tau}(x)$ be the 
generating function of the sequence $a_{\tau}(n)$, that is,  
		$$a_{\tau}(x)=\sum_{n\geq 0} a_{\tau}(n)x^n.$$ 
We find an explicit expression for the generating function $a_{\tau}(x)$.

\begin{theorem}
\label{a12}
Let $\tau\in S_k(123,132)$. Then:
\begin{enumerate}
\item [(i)] 	there exist $r_1,r_2,\dots,r_m\geq 1$ with $r_1+r_2+\dots+r_m=k$ 
		such that $$\tau=(\beta_1, \beta_2 , \dots, \beta_m),$$
		where $\beta_i=(t_i-1, t_i-2,\dots, t_i-r_i+1, t_i)$, and 
		$t_i=k-(r_1+\dots+r_{i-1})$ for $i=1,2,\dots,m$;

\item[(ii)] 	$$a_{\tau}(x)=
		\left\bracevert
		\begin{array}{ccccc}
		f_{r_1}(x)	&	-g_{r_1}(x)	&	0	&	\dots	&	0 \\
		f_{r_2}(x)	&	1		&   -g_{r_2}(x)	&	\ddots	&	0 \\
		f_{r_3}(x)	&	0		&	1	&  	\ddots	& \vdots \\
		\vdots		&	\vdots		&	\ddots	&	\ddots	&	-g_{r_{m-1}}(x)	\\
		f_{r_m}(x)	&	0		&	0	&	0	&	1 
		\end{array}
		\right\bracevert,$$
where $f_d(x)=\frac{1-x}{1-2x+x^d}$ and $g_d(x)=\frac{x^d}{1-2x+x^d}$.
\end{enumerate}
\end{theorem}
\begin{proof}

$(i)$ Let $\tau\in S_k(123,132)$; choose $r_1$ such that $\tau_{r_1}=k$. 
Since $\tau$ avoids $132$, we see that $\tau_j\geq k-r_1+1$ for all 
$j\leq r_1$, and since $\tau$ avoids $123$, we get that 
$\tau=(\beta_1,\tau')$, where $\tau'\in S_{k-r_1}(123,132)$, and so on.\\

$(ii)$ Let $\tau=(\beta_1,\dots,\beta_m)$, and let $\alpha\in S_n(123,132,\tau)$;  
choose $t$ such that $\alpha_t=n$. Similarly to $(i)$,
$\alpha=(n-1,\dots,n-t+1,n,\alpha_{t+1},\dots,\alpha_n)$, therefore 
$$a_{\tau}(n)=\sum_{t=1}^{r_1-1} a_{\tau}(n-t) +\sum_{t=r_1}^n a_{\tau'}(n-t),$$
which means
	$$a_{\tau}(n)=2a_{\tau}(n-1)-a_{\tau}(n-r_1)+a_{\tau'}(n-r_1)$$
for all $n\geq k+1$, where $\tau'=(\beta_2,\dots,\beta_m)$. Hence 
	$$a_{(\beta_i,\dots,\beta_m)}(n)=2a_{(\beta_i,\dots,\beta_m)}(n-1)-a_{(\beta_i,\dots,\beta_m)}(n-r_i)+a_{(\beta_{i+1},\dots,\beta_m)}(n-r_i)$$
for all $n\geq t_i+1$, or equivalently,
	$$\begin{array}{ll}
	\sum\limits_{n\geq t_i+1} a_{(\beta_i,\dots,\beta_m)}(n)x^n
		&=2x \sum\limits_{n\geq t_i} a_{(\beta_i,\dots,\beta_m)}(n)x^n \\
		&-x^{r_i}\sum\limits_{n\geq t_i+1} a_{(\beta_i,\dots,\beta_m)}(n)x^n \\
		&+x^{r_i}\sum\limits_{n\geq t_i+1} a_{(\beta_{i+1},\dots,\beta_m)}(n)x^n.
	\end{array}$$
Since $a_{(\beta_i,\dots,\beta_m)}(n)=2^{n-1}$ for all $n\leq t_i-1$, 
$a_{(\beta+i,\dots,\beta_m)}(t_i)=2^{t_i-1}-1$, and 
$a_{(\beta_i,\dots,\beta_m)}(0)=1$ (Simion and Schmidt ~\cite{simion}, Proposition $7$), 
we obtain that  
	$$(1-2x+x^{r_i})a_{(\beta_i,\dots,\beta_m)}(x)=1-x+x^{r_i}a_{(\beta_{i+1},\dots,\beta_m)}(x)$$
for all $i\leq m-1$, and 
		$$(1-2x+x^{r_m})a_{(\beta_m)}(x)=1-x.$$
Hence, by Lemma \ref{det}, the theorem holds. \hfill\rule{2mm}{2mm}
\end{proof}

\begin{example}
\label{ex11}
Let $T=\{123,132\}$. By Theorem \ref{a12},
\begin{enumerate}
\item 	$|S_n(T,3214)|=t_n$, where $t_n$ is the $n$-th Tribonacci 
	number ~\cite{feinberg}, and 
 	$|S_n(T,3241)|=f_{n+2}-1$, where $f_n$ is the $n$-th 
	Fibonacci number. \\

\item 	$|S_n(T,3412)|=|S_n(T,4231)|={n\choose 2}+1$.\\

\item 	$|S_n(T,3421)|=3n-5$.
\end{enumerate}
\end{example}

\subsection{$T=\{123,\ 231\}$.}

In this subsection, we calculate the cardinality of the set 
$S_n(123,231,\tau)$, where $\tau\in S_k(123,231)$. This cardinality 
we denote by $b_{\tau}(n)$.

\begin{lemma}
\label{l14}
Let $\tau\in S_k(123,231)$. Then, either there exists $r$, 
$1\leq r\leq k-1$, such that $\tau=(r,\dots,2,1,k,k-1,\dots,r+1)$, 
and hence
   	$$b_{\tau}(n)=(k-2)n-\frac{k(k-3)}{2}\ \ \mbox{for all}\ n\geq k,$$
or $\tau=(k,\tau')\neq (k,\dots,2,1)$ such that $\tau'\in S_{k-1}(123,231)$,
and hence
	$$b_{\tau}(n)=b_{\tau'}(n-1)+n-1\ \ \mbox{for all}\ n\geq k.$$
\end{lemma}
\begin{proof}

Let $\tau\in S_k(123,231)$; put $r=\tau_1$. Since $\tau$ avoids 
$123$, we see that $\tau$ contains $(r,k,\dots,r+1)$, since $\tau$ 
avoids $231$, we get that $\tau=(r,\tau',k,k-1,\dots,r+1)$, and since 
$\tau$ avoids $123$, we have two cases: either $\tau=(r,\dots,1,k,\dots,r+1)$ 
for $1\leq r\leq k-1$, or $\tau=(k,\tau')$ such that 
$\tau'\in S_{k-1}(123,231)$. Now let us consider the two cases:
\begin{enumerate}
\item Let $\alpha\in S_n(123,231,\tau)$, where $\tau=(r,\dots,1,k,\dots,r+1)$, 
$1\leq r\leq k-1$. Similarly to the above, we have two cases for $\alpha$: in the 
first case $\alpha=(t,\dots,1,n,n-1,\dots,t+1)$ for $1\leq t\leq n-1$, so there are 
$k-2$ permutations like $\alpha$. In the second case 
$\alpha=(n,\alpha_2,\dots,\alpha_n)$, so there are $b_{\tau}(n-1)$ 
permutations, which means 
$b_{\tau}(n)=b_{\tau}(n-1)+k-2$. Besides, 
$a_{\tau}(k)=k(k-1)/2$ (see Simion and Schmidt ~\cite{simion}, Proposition $11$), 
hence $b_{\tau}(n)=(k-2)n-\frac{k(k-3)}{2}$.

\item	Let $\alpha\in S_n(123,231,\tau)$, and let $\tau=(k,\tau')\neq (k,\dots,1)$ such that  
$\tau'\in S_{k-1}(123,231)$. Similarly to the above, we have two cases for 
$\alpha$: in the first case $\alpha=(t,\dots,1,n,n-1,\dots,t+1)$ for 
$1\leq t\leq n-1$, so there are $n-1$ permutations like $\alpha$. In the 
second case $\alpha=(n,\alpha_2,\dots,\alpha_n)$, so there are 
$b_{\tau'}(n-1)$ permutations. Hence $b_{\tau}(n)=b_{\tau'}(n-1)+n-1$.
\end{enumerate}
\hfill\rule{2mm}{2mm}
\end{proof}

\begin{theorem}
\label{t14}
Let $\tau\in S_k(123,231)$. Then :
\begin{enumerate}
\item[(i)]	there exist $m$, $2\leq m\leq k+1$, and $r$, $1\leq r\leq m-2$, 
		such that 
		$$\tau=(k,\dots,m,r,\dots,1,m-1,\dots,r+1);$$

\item[(ii)]	for all $n\geq k$
		  $$b_{\tau}(n)=(k-2)n - \frac{k(k-3)}2.$$
\end{enumerate}
\end{theorem}
\begin{proof}

$(i)$ Immediately, by Lemma \ref{l14},
	$$\tau=(k,\dots,m,r,\dots,1,m-1,\dots,r+1),$$ 
where $2\leq m\leq k+1$, $1\leq r\leq m-2$. 

$(ii)$ Again, by Lemma \ref{l14}, for all $n\geq k$
  $$b_{\tau}(n)=\sum_{j=1}^{k-m+1} (n-j) + (m-3)(n-(k-m+1))+\frac{(m-1)(4-m)}{2},$$
hence, this theorem holds. \hfill\rule{2mm}{2mm}
\end{proof}

\begin{example}
\label{ex12}
Let $T=\{123,231\}$; by the Theorem \ref{t14}  
 $$|S_n(T,4312)|=|S_n(T,1432)|=|S_n(T,2143)|=|S_n(T,3214)|=2n-2.$$
\end{example}
\subsection{$T=\{132,\ 213\}$.}

Let $c_{\tau}(n)=|S_n(132,213,\tau)|$, and let $c_{\tau}(x)$ be the 
generating function of the sequence $c_{\tau}(n)$. 
We find an explicit expression for the generating function $c_{\tau}(x)$.

\begin{theorem}
\label{t23}
Let $\tau\in S_k(132,213)$. Then:
\begin{enumerate}

\item[(i)] 	there exist $k+1=r_0>r_1>\dots>r_m\geq 1$ such that 
		$$\tau=(r_1,r_1+1,\dots,k,r_2,r_2+1,\dots,r_1-1,\dots,r_m,r_m+1,\dots,r_{m-1}-1);$$

\item[(ii)] 	$$c_{\tau}(x)=
		\left\bracevert
		\begin{array}{ccccc}
		f_{r_0-r_1}(x)	&    -g_{r_0-r_1}(x)	&	0	&	\dots	&	0 \\
		f_{r_1-r_2}(x)	&	1		&   -g_{r_1-r_2}(x)&	\ddots	&	0 \\
		f_{r_2-r_3}(x)	&	0		&	1	&  	\ddots	& \vdots \\
		\vdots		&	\vdots		&	\ddots	&	\ddots	&	-g_{r_{m-2}-r_{m-1}}(x)	\\
		f_{r_{m-1}-r_m}(x)&	0		&	0	&	0	&	1 
		\end{array}
		\right\bracevert,$$
where $f_d(x)$ and $g_d(x)$ have the same meaning as in Theorem \ref{a12}.
\end{enumerate}
\end{theorem}
\begin{proof}

$(i)$ Let $\tau\in S_k(132,213)$, and let $r_1=\tau_1$. 
Since $\tau$ avoids $132$, we see that $\tau$ contains $(r_1,r_1+1,\dots,k)$, 
and since $\tau$ avoids $213$, we get that $\tau=(r_1,r_1+1,\dots,k,\tau')$, 
where $\tau'\in S_{r_1-1}(132,213)$, and so on.\\

$(ii)$ Let $\alpha\in S_n(132,213,\tau)$, and let $t=\alpha_1$; similarly to $(i)$,
$\alpha=(t,t+1,\dots,n,\alpha_{n-t+2},\dots,\alpha_n)$.
Therefore, for $t\leq n-k+r_1$ we have 
$\alpha\in S_n(132,213,\tau)$ if and only if 
$(\alpha_{n-t+2},\dots,\alpha_n)\in S_{t-1}(132,213,\tau')$, and for 
$t\geq n-k+r_1+1$ we have $\alpha\in S_n(132,213,\tau)$ if and 
only if $(\alpha_{n-t+2},\dots,\alpha_n)\in S_{t-1}(132,213,\tau)$. Hence 
$$c_{\tau}(n)=\sum_{t=1}^{n-k+r_1} c_{\tau'}(t-1) + \sum_{t=n-k+r_1+1}^n c_{\tau}(t-1),$$
which means that
$$c_{\tau}(n)=2c_{\tau}(n-1) -c_{\tau}(n-k+r_1-1)+c_{\tau'}(n-k+r_1-1).$$
Let us define $t_i=r_{i-1}-r_i$ for $i=1,2,\dots,m$, so 
	$$c_{\tau}(n)=2c_{\tau}(n-1)-c_{\tau}(n-t_1)+c_{\tau'}(n-t_1).$$
If $\tau=(1,2,\dots,k)$, then immediately $c_{\tau}(x)=\frac{1-x}{1-2x+x^k}$, 
hence, by Lemma \ref{det}, this theorem holds. \hfill\rule{2mm}{2mm}
\end{proof}

\begin{corollary}
\label{cc2}
Let $k\geq 2$. For all $n\geq 0$, 
	$$c_{(k,\dots,2,1)}(n)=\sum_{j=0}^{k-2} {{n-1}\choose j}.$$
\end{corollary}
\begin{proof}
By the proof of Theorem \ref{t23},
 $$c_{(k,\dots,2,1)}(n)=2^k + \sum_{t=k-1}^{n-1} c_{(k-1,\dots,2,1)}(t),$$
which means that 
   $$c_{(k,\dots,2,1)}(n)=c_{(k,\dots,2,1)}(n-1)+c_{(k-1,\dots,2,1)}(n-1).$$
Besides, $c_{(k,\dots,2,1)}(k)=2^{k-1}-1$, and 
$c_{(k,\dots,2,1)}(n)=2^{n-1}$ for $1\leq n\leq k-1$
(see Simion and Schmidt ~\cite{simion}, Proposition $8$). 
Hence, the corollary is true. \hfill\rule{2mm}{2mm}
\end{proof}

\begin{example}
\label{ex13}
Let $T=\{132,213\}$.
\begin{enumerate}
\item 	By Corollary \ref{cc2}, $|S_n(T,4321)|={n\choose 2}+1$. \\

\item	By Theorem \ref{t23}, $|S_n(T,1234)|=t_n$ where $t_n$ is the 
	$n$-th Tribonacci number. \\

\item	By Theorem \ref{t23}, $|S_n(T,2341)|=f_{n+2}-1$ where $f_n$ is 
	the $n$-th Fibonacci number.\\

\item  By Theorem \ref{t23}, $|S_n(T,3412)|=|S_n(T,3421)|=|S_n(T,4231)|={n\choose 2}+1$.
\end{enumerate}
\end{example}
\subsection{$T=\{213,\ 231\}$.}

Let $d_{\tau}(n)=|S_n(213,231,\tau)|$, and let $d_{\tau}(x)$ be the 
generating function of the sequence $d_{\tau}(n)$. We find an 
explicit expression for the generating function $d_{\tau}(x)$.

\begin{theorem}
\label{t24}
Let $\tau\in S_k(213,231)$. Then:
\begin{enumerate}
\item[(i)]	$\tau_i$ is either the right maximum, or
		the righr minimum, for all $1\leq i\leq k-1$;

\item[(ii)] 	$$d_{\tau}(x)=
		\left\bracevert
		\begin{array}{ccccc}
		1	&    -g(x)	&	0	&	\dots	&	0 \\
		1	&	1	&   -g(x)	&	\ddots	&	0 \\
		1	&	0	&	1	&  	\ddots	& \vdots \\
		\vdots	&	\vdots	&	\ddots	&	\ddots	&	0	\\
		1	&	0	&	0	&	0	&	-g(x) \\
		1	&	0	&	0	&	0	&	1	
		\end{array}
		\right\bracevert,$$
where $g(x)=\frac{x}{1-x}$.
\end{enumerate}
\end{theorem}
\begin{proof}

$(i)$ Let $\tau\in S_k(213,231)$; if $2\leq \tau_1\leq k-1$, then 
$\tau$ contains either $(\tau_1,1,k)$ or $(\tau_1,k,1)$, which means 
$\tau$ contains either $213$ or $231$, hence $\tau_1=1$ or $\tau_1=k$, 
and so on. \\

$(ii)$ Let $\alpha\in S_n(213,231,\tau)$; similarly to $(i)$, 
$\alpha_1=1$ or $\alpha_1=n$. Let $\tau=(\tau_1,\tau')$, hence 
in the above two cases ($\tau_1=1$ or $\tau_1=k$) we obtain 
	$$d_{\tau}(n)=d_{\tau}(n-1)+d_{\tau'}(n-1)$$
for all $n\geq k$. Besides, $d_{\tau}(0)=1$, 
$d_{\tau}(k)=2^{k-1}-1$, and $d_{\tau}(n)=2^{n-1}$ for all $1\leq n\leq k-1$ 
(see Simion and Schmidt ~\cite{simion}, Proposition $10$). Hence, 
similarly to Theorem \ref{a12}, the theorem holds.
\hfill\rule{2mm}{2mm}
\end{proof}

Immediately by Theorem \ref{t24},
	$$|S_n(123,132,231)|=|S_n(132,231,321)|=n,$$
for all $n\geq 0$, which means, we have a generalization 
Proposition $16$ and Lemma $6(b)$ of Simion and Schmidt 
~\cite{simion}. 

\begin{example}
\label{ex14}
Let $T=\{213,231\}$. By Theorem \ref{t24},
  $$|S_n(T,1234)|=|S_n(T,1243)|=|S_n(T,1423)|=|S_n(T,1432)|={n\choose 2}+1.$$ 
\end{example}
\section{Three patterns from $S_3$ and a pattern from $S_k$}

In this section, we calculate the cardinality of the sets 
$S_n(T,\tau)$ such that $T\subset S_3$, $|T|=3$ and $\tau\in S_k(T)$ 
for $k\geq 3$. By Remark \ref{rem} and by three natural operations
the complementation, the reversal and the inverse 
(see Simion and Schmidt ~\cite{simion}, Lemma $1$), 
we have to consider the following five possibilities:\\

$\begin{array}{rlll}
  1)&	S_n(123,132,213,\tau),	& \ \mbox{where}	& \ \tau\in S_k(123,132,213), \\
  2)&	S_n(123,132,231,\tau), 	& \ \mbox{where}	& \ \tau\in S_k(123,132,231), \\
  3)&	S_n(123,213,231,\tau),	& \ \mbox{where}	& \ \tau\in S_k(123,213,231), \\
  4)&	S_n(123,231,312,\tau),	& \ \mbox{where} 	& \ \tau\in S_k(123,231,312), \\
  5)&	S_n(132,213,231,\tau),	& \ \mbox{where} 	& \ \tau\in S_k(132,213,231). 
\end{array}$
\\
\begin{remark}
By Erd\"os and Szekeres ~\cite{erdos}, $|S_n((1,2,\dots,a),(b,b-1,\dots,1))|=0$ 
for all $n\geq (a-1)(b-1)+1$, where $a,b\geq 1$.  
Therefore, in what follows we assume that $\tau\in S_k(T)$ and 
$\tau\neq (k,k-1,\dots,1)$, since $123\in T$.
\end{remark}

The main body of this section is divided into five subsections corresponding 
to the above five cases.
\subsection{$T=\{123,\ 132,\ 213\}$.}
Let $e_{\tau}(x)$ be the generating function of the sequence 
$|S_n(T,\tau)|$.  We find  an explicit expression for the 
generating function $e_{\tau}(x)$.
\begin{lemma}
\label{l123}
Let $T=\{123,132,213\}$ and $\tau\in S_k(T)$. Then, either 
there exists $\tau'\in S_{k-1}(T)$ such that $\tau=(k,\tau')\neq (k,k-1,\dots,1)$, 
and hence
	 $$|S_n(T,\tau)|=|S_{n-1}(T,\tau')|+|S_{n-2}(T,\tau')|\ \ \mbox{for any}\ n\geq k,$$
or there exists $\tau''\in S_{k-2}(T)$ such that $\tau=(k-1,k,\tau')$, and
hence
 	$$|S_n(T,\tau)|=|S_{n-1}(T,\tau)|+|S_{n-2}(T,\tau'')|\ \ \mbox{for any}\ n\geq k.$$
\end{lemma}
\begin{proof}

Let $\tau\in S_k(T)$; since $\tau$ avoids $123$ and $132$ we have 
either $\tau_1=k$ or $\tau_1=k-1$. If $\tau_1=k-1$, then, 
since $\tau$ avoids $213$, we see that $\tau=(k-1,k,\tau'')$.
Now we consider the two cases:
\begin{enumerate}
\item	Let $\tau=(k,\tau')$, $\alpha\in S_n(T,\tau)$. Similarly to the above, 
either $\alpha=(n,\alpha_2,\dots,\alpha_n)$, or 
$\alpha=(n-1,n,\alpha_3,\dots,\alpha_n)$, so evidently 
	 $$|S_n(T,\tau)|=|S_{n-1}(T,\tau')|+|S_{n-2}(T,\tau')|.$$
	
\item Let $\tau=(k-1,k,\tau'')$, $\alpha\in S_n(T,\tau)$. Similarly to the above,
either $\alpha=(n,\alpha_2,\dots,\alpha_n)$, or 
$\alpha=(n-1,n,\alpha_3,\dots,\alpha_n)$, so evidently 
	 $$|S_n(T,\tau)|=|S_{n-1}(T,\tau)|+|S_{n-2}(T,\tau'')|.$$
\end{enumerate}
\hfill\rule{2mm}{2mm}
\end{proof}

For any permutation $\tau\in S_k$ such that $\tau_1=k$ we define  
$p(\tau)=1$ and $q(\tau)=(\tau_2,\dots,\tau_k)$, and for 
any permutation $\tau\in S_k$ such that $\tau_1=k-1$, and $\tau_2=k$ 
we define $p(\tau)=2$ and $q(\tau)=(\tau_3,\dots,\tau_k)$.
Also, let $m(\tau)=(m_1,m_2,\dots,m_r)$ where $m_i=p(q^{i-1}(\tau)))$ 
for $1\leq i\leq r$, 
and $q^0(\tau)=\tau$, $q^{i-1}(\tau)=q(q^{i-2}(\tau)$ for 
$i\geq 2$.

\begin{theorem}
\label{t123}
Let $T=\{123,132,213\}$, $\tau\in S_k(T)$, and $m(\tau)=(m_1,\dots,m_r)$.
Then 
	$$e_{\tau}(x)=
		\left\bracevert
		\begin{array}{ccccc}
		u_{m_1}(x)	&    -v_{m_1}(x)	&	0	&	\dots	&	0 \\
		u_{m_2}(x)	&	1		&   -v_{m_2}(x) &	\ddots	&	0 \\
		u_{m_3}(x)	&	0		&	1	&  	\ddots	& \vdots \\
		\vdots		&	\vdots		&	\ddots	&	\ddots	&	-v_{m_{r-1}}(x)	\\
		u_{m_r}(x)&	0		&	0	&	0	&	1 
		\end{array}
		\right\bracevert,$$
where $u_1(x)=1$, $u_2(x)=\frac1{1-x}$, $v_1(x)=x(1+x)$, 
and $v_2(x)=\frac{x^2}{1-x}$.	
\end{theorem}
\begin{proof}
Let $\tau\in S_k(T)$; by Lemma \ref{l123}, there are two cases:
\begin{enumerate}
\item	$\tau_1=k$. So $p(\tau)=1$, $q(\tau)=(\tau_2,\dots,\tau_k)$,  
	and for all $n\geq k$
	 $$|S_n(T,\tau)|=|S_{n-1}(T,\tau^2)|+|S_{n-2}(T,\tau^2)|.$$
	Besides, $|S_n(T,\tau)|=f_n$ for all $n\leq k-1$ 
	(see Simion and Schmidt ~\cite{simion}, Proposition $15$), where 
	$f_n$ is the $n$-th Fibonacci number. Hence 
	$$e_{\tau}(x)=x(1+x)e_{q(\tau)}(x)+1.$$

\item 	$\tau_1=k-1$. So $p(\tau)=2$, $q(\tau)=(\tau_3,\dots,\tau_k)$,
	and for all $n\geq k$, 
	 $$|S_n(T,\tau)|=|S_{n-1}(T,\tau)|+|S_{n-2}(T,\tau'')|.$$
	Besides, $|S_n(T,\tau)|=f_n$ for all $n\leq k-1$ 
	(see Simion and Schmidt ~\cite{simion}, Proposition $15$), where 
	$f_n$ is the $n$-th Fibonacci number. Hence 
	$$e_{\tau}(x)=\frac{x^2}{1-x}e_{q(\tau)}(x)+\frac{1}{1-x}.$$	
\end{enumerate}
Hence, by the definitions and Lemma \ref{det}, the theorem holds.
\hfill\rule{2mm}{2mm}
\end{proof}

\begin{example}
\label{ex21}
Let $T=\{123,132,213\}$. By Theorem \ref{t123},
\begin{enumerate}
\item 	$|S_n(T,3412)|=n$.

\item	$|S_n(T,4231)|=|S_n(T,3421)|=4$.
\end{enumerate}
\end{example}
\subsection{T=\{123,\ 132,\ 231\}.}
\begin{theorem}
\label{t124}
Let $T=\{123,132,231\}$ and $\tau\in S_k(T)$. Then: 
\begin{enumerate}
\item[(i)]	there exists $r$, $1\leq r\leq k$, 
	such that $\tau=(k,\dots,r+1,r-1,\dots,1,r)$; 

\item[(ii)]	for all $n\geq k$
		$$|S_n(T,(k,\dots,r+1,r-1,\dots,1,r)|=k-1,$$
		where $2\leq r\leq k$.
\end{enumerate}
\end{theorem}
\begin{proof}

$(i)$ Let $\tau\in S_k(T)$; put $r=\tau_n$. Since $\tau$ avoids $123$, 
we see that $\tau$ contains $(r-1,\dots,1,r)$, since $\tau$ avoids 
$132$, we see that $\tau=(\tau_1,\dots,\tau_{k-r},r-1,\dots,1,r)$, 
and since $\tau$ avoids $231$, we get that 
$\tau=(k,\dots,r+1,r-1,\dots,1,r)$.\\

$(ii)$ Let $\alpha\in S_n(T,\tau)$; similarly to $(i)$, 
$\alpha=(n,\dots,t+1,t-1,\dots,1,t))$ for $1\leq t\leq n$, 
hence $|S_n(T,\tau)|=k-1.$ \hfill\rule{2mm}{2mm}
\end{proof}

\begin{example}
\label{ex22}
Let $T=\{123,132,231\}$; by Theorem \ref{t124},
	$$|S_n(T,4312)|=|S_n(T,4213)|=|S_n(T,3214)|=3.$$
\end{example}
\subsection{T=\{123,\ 213,\ 231\}.}
\begin{theorem}
\label{t134}
Let $T=\{123,213,231\}$ and $\tau\in S_k(T)$. Then: 
\begin{enumerate}
\item[(i)]	there exists $r$, $1\leq r\leq k$, 
	such that $\tau=(k,\dots,r+1,1,r,\dots,2)$; 

\item[(ii)] for all $n\geq k$
      $$|S_n(T,(k,\dots,r+1,1,r,\dots,2))|=k-1,$$
	where $2\leq r\leq k$.
\end{enumerate}
\end{theorem}
\begin{proof}

$(i)$ Let $\tau\in S_k(T)$ and choose $r$ such that $\alpha_{k-r+1}=1$. 
Since $\tau$ avoids $213$, we get that $\tau_i>\tau_j$ for all $i<k-r+1<j$, 
since $\tau$ avoids, $123$ we see that $\tau$ contains $(1,r,\dots,2)$, 
and since $\alpha$ avoids $231$, we get that $\tau=(k,\dots,r+1,1,r,\dots,2)$.\\

$(ii)$ Let $\alpha\in S_n(T,\tau)$; similarly to $(i)$,
$\alpha=(n,\dots,t+1,1,t,\dots,2)$ for $1\leq t\leq n$, 
hence $|S_n(T,\tau)|=k-1$.\hfill\rule{2mm}{2mm}
\end{proof}

\begin{example}
\label{ex23}
Let $T=\{123,213,231\}$; by Theorem \ref{t134},
	$$|S_n(T,4312)|=|S_n(T,4132)|=|S_n(T,1432)|=3.$$
\end{example}
\subsection{T=\{123,\ 231,\ 312\}.}
\begin{theorem}
\label{t145}
Let $T=\{123,231,312\}$ and $\tau\in S_k(T)$. Then  :
\begin{enumerate}
\item[(i)]	there exists $r$, $1\leq r\leq k$, 
	such that $\tau=(r,\dots,2,1,k,\dots,r+1)$; 

\item[(ii)]	for all $n\geq k$
	$$|S_n(T,(r,\dots,2,1,k,\dots,r+1))|=k-1,$$
	where $1\leq r\leq k-1$.
\end{enumerate}
\end{theorem}
\begin{proof}

$(i)$ Let $\tau\in S_k(T)$; put $r=\tau_1$. Since $\tau$ avoids $123$, we 
get that $\tau$ contains $(r,k,\dots,r+1)$, and since $\tau$ avoids 
$231$, we see that $\tau=(r,\dots,k,\dots,r+2,r+1)$, and since $\tau$ 
avoids $312$, we get that $\tau=(r,\dots,2,1,k,\dots,r+2,r+1)$ for 
$r=1,\dots,k$.\\

$(ii)$ Let $r\leq k-1$ and $\alpha\in S_n(T,\tau)$; similarly to $(i)$,
$\alpha=(t,\dots,2,1,n,\dots,t+2,t+1)$ for $1\leq t\leq n$. 
Hence $|S_n(T,\tau)|=k-1$. \hfill\rule{2mm}{2mm}
\end{proof}

\begin{example}
\label{ex24}
Let $T=\{123,231,312\}$; by Theorem \ref{t145},
	$$|S_n(T,1432)|=|S_n(T,2143)|=|S_n(T,3214)|=3.$$
\end{example}
\subsection{T=\{132,\ 213,\ 231\}.}
\begin{theorem}
\label{t234}
Let $T=\{132,213,231\}$ and $\tau\in S_k(T)$. Then :
\begin{enumerate}
\item[(i)]	there exists $r$, $1\leq r\leq k$ ,
	such that $\tau=(k,\dots,r+1,1,2,\dots,r)$; 

\item[(ii)]	for all $n\geq k$
	 $$|S_n(T,(k,\dots,r+1,1,2,\dots,r))|=k-1.$$
\end{enumerate}
\end{theorem}
\begin{proof}

$(i)$ Let $\tau\in S_k(T)$; put $r=\tau_n$. Since $\tau$ avoids $231$, we 
get that $\tau$ contains $(k,\dots,r+1,r)$, since $\tau$ avoids 
$132$, we see that $\tau=(k,\dots,r+1,\tau_{k-r+1},\dots,\tau_{k-1},r)$,  
and since $\tau$ avoids $213$, we get that 
$\tau=(k,k-1,\dots,r+1,1,2,\dots,r)$ for $r=1,\dots,k$.\\

$(ii)$ Let $\alpha\in S_n(T,\tau)$; similarly to (i), 
$\alpha=(n,\dots,t+1,1,2,\dots,t)$ for $1\leq t\leq n$. 
Hence $|S_n(T,\tau)|=k-1$. \hfill\rule{2mm}{2mm}
\end{proof}

\begin{example}
\label{ex25}
Let $T=\{132,213,231\}$; by Theorem \ref{t234},
    $$|S_n(T,4321)|=|S_n(T,4312)|=|S_n(T,4123)|=|S_n(T,1234)|=3.$$
\end{example}

{\begin{table}[p]
    \begin{tabular}{|l|l|l|l|} \hline
					&		&							&			\\
	\emph{Wilf class $C$} 		& \emph{$|C|$}	& \emph{Cardinality of $S_n(T)$, $T\in C$}		& \emph{Reference}	\\ \hline\hline

					&		&							&			\\
	$\overline{\{123,4321\}}$	& $490$		& $0$							& Erd\"os and Szekeres ~\cite{erdos}	\\ \hline
					&		&							&			\\
	$\overline{\{123,1234\}}$	& $60$ 		& $c_n$							& West ~\cite{west}, Knuth ~\cite{knuth} \\ \hline
					&		&							&			\\
	$\overline{\{123,1432\}}$	& $46$		& $f_{2n-2}$ 						& West ~\cite{west}	\\ \hline
					&		&							&			\\
	$\overline{\{132,3421\}}$	& $12$ 		& $1+(n-1)2^{n-2}$					& West ~\cite{west} , Guibert ~\cite{guibert}	\\ \hline
					&		&							&			\\
	$\overline{\{123,2431\}}$	& $8$		& $3\cdot 2^{n-1}-{{n+1}\choose{2}}-1$			& West ~\cite{west}	\\ \hline
					&		&							&			\\
	$\overline{\{123,3421\}}$	& $4$		& ${n\choose 4}+2{n\choose 3}+n$			& West ~\cite{west}	\\ \hline
					&		&							&			\\
	$\overline{\{132,3214\}}$	& $4$		& Generating function $\frac{(1-x)^3}{1-4x+5x^2-3x^3}$	& West ~\cite{west}	\\ \hline
					&		&							&			\\
	$\overline{\{132,4321\}}$	& $4$		& ${n\choose 4}+{{n+1}\choose 4}+{n\choose 2}+1$	& West ~\cite{west}	\\ \hline 
					&		&							&			\\
	$\overline{\{123,3412\}}$	& $2$		& $2^{n+1}-{{n+1}\choose {3}}-2n-1$			& Billey, Jockusch and Stanley ~\cite{billey}	\\ \hline
					&		&							&			\\
	$\overline{\{123,4231\}}$	& $2$		& ${n\choose 5}+2{n\choose 4}+{n\choose 3}+{n\choose 2}+1$ & West ~\cite{west}	\\ \hline\hline
 					&		&							&			\\
	$\overline{\{123,132,1234\}}$ 	& $160$ 	& $2^{n-1}$						& West ~\cite{west}, Simion and Schmidt ~\cite{simion}	\\ \hline 
					&		&							&			\\
	$\overline{\{123,132,3412\}}$	& $118$ 	& ${n\choose 2}+1$					& Section $2$, Examples \ref{ex11}, \ref{ex13}, \ref{ex14}	\\ \hline
					&		&							&			\\
	$\overline{\{123,312,1432\}}$ 	& $24$		& $2n-2$						& Section $2$, Example \ref{ex12}  		\\ \hline
					&		&							&			\\
	$\overline{\{123,132,3241\}}$	& $12$ 		& $f_{n+2}-1$						& Section $2$, Examples \ref{ex11}, \ref{ex13}	\\ \hline
					&		&							&			\\
	$\overline{\{123,132,3421\}}$	& $8$  		& $3n-5$						& Section $2$, Example \ref{ex11}		\\ \hline
					&		&							&			\\
	$\overline{\{123,132,3214\}}$	& $6$  		& $t_n$ 						& Section $2$, Examples \ref{ex11}, \ref{ex13}	\\ \hline\hline
						&		&						&			\\
	$\overline{\{123,132,231,1234\}}$	& $282$		& $n$						& West ~\cite{west}, Simion and Schmidt ~\cite{simion}  \\ \hline 
						&		&						&			\\
	$\overline{\{123,132,231,3214\}}$	& $46$		& $3$						& Section $3$, Examples \ref{ex22}, \ref{ex23}, \ref{ex24}, \ref{ex25}	\\ \hline
						&		&						&			\\
	$\overline{\{123,132,213,1234\}}$	& $38$ 		& $f_{n+1}$					& West ~\cite{west}, Simion and Schmidt ~\cite{simion} \\ \hline 
						&		&						&			\\
	$\overline{\{123,132,213,3421\}}$	& $6$		& $4$						& Section $3$, Example \ref{ex21}	\\ \hline\hline
						&		&						&			\\
	$\overline{\{123,132,213,231,1234\}}$	& $100$ 	& $2$						& Section $4$, Theorem \ref{t1234}	\\ \hline
						&		&						&			\\
	$\overline{\{123,132,213,231,4312\}}$ 	& $56$		& $1$						& Section $4$, Theorem \ref{t1234}	\\ \hline

    \end{tabular}
\caption {Wilf classes of $\{T,\tau\}$, where $T\subseteq S_3$, $\tau\in S_4$}
\label{tab34}
\end{table} }

\section{At least four patterns from $S_3$ and a pattern from $S_k$}

By Simion and Schmidt ~\cite{simion}, Proposition $17$, 
		$$|S_n(T)|=2,$$ 
where $\{123,321\}\not\subset T\subset S_3$, $|T|=4,5$, and 
		$$|S_n(T)|=0,$$ 
where $\{123,321\}\subset T$. Hence, we obtain the following theorem.

\begin{theorem}
\label{t1234}
Let $T\subset S_3$, $|T|\geq 4$ and $\tau\in S_k(T)$. For all $n\geq k$
$$|S_n(T,\tau)|=\left\{ 
	\begin{array}{ll}
	2-\delta_{\tau,(1,2,\dots,k)} -\delta_{\tau,(k,\dots,2,1)}	& 123,321\notin T \\
	2-\delta_{\tau,(k,\dots,2,1)}					& 123\in T,\ 321\notin T \\
	2-\delta_{\tau,(1,2,\dots,k)}					& 123\notin T,\ 321\in T \\
	0								& 123,321\in T.  
	\end{array}
	\right.,$$
where $\delta_{x,y}$ the Kronecker symbol.
\end{theorem}
\section{Wilf classes of $\{T,\tau\}$, where $T\subseteq S_3$, $\tau\in S_4$}

By all the examples in all the sections we obtain Table \ref{tab34}. 
This table describes all the Wilf classes of sets of permutations avoiding a 
pattern from $S_4$ and a set of patterns from $S_3$. It contains $22$ 
Wilf classes for sets $\{T,\tau\}$ where $T\subseteq S_3$, 
$\tau\in S_4$.

\end{document}